\documentclass[journal,twoside]{IEEEtran}
\usepackage{cite}

\ifCLASSINFOpdf
\else
\fi
\usepackage[cmex10]{amsmath}
\usepackage{array}

\usepackage{mdwmath}
\usepackage{mdwtab}
\usepackage[caption=false,font=footnotesize]{subfig}
\usepackage{fixltx2e}

\usepackage{stfloats}
\usepackage{url}

\usepackage[american,fulldiode]{circuitikz} 
\usepackage{graphicx} 
\usepackage{multirow}
\usepackage{placeins}
\usepackage{amsfonts}

\begin{document}
%
\title{A Sufficient Condition for Power Flow Insolvability with Applications to Voltage Stability Margins}
%
%
%

\author{Daniel K. Molzahn, \IEEEmembership{Student Member, IEEE},  Bernard C. Lesieutre, \IEEEmembership{Senior Member, IEEE}, \\ and Christopher L. DeMarco, \IEEEmembership{Member, IEEE}
\thanks{University of Wisconsin-Madison Department of Electrical and Computer Engineering: molzahn@wisc.edu, lesieutre@engr.wisc.edu, demarco@engr.wisc.edu}}

\maketitle

\begin{abstract}
For the nonlinear power flow problem specified with standard PQ, PV, and slack bus equality constraints, we present a sufficient condition under which the specified set of nonlinear algebraic equations has no solution. This sufficient condition is constructed in a framework of an associated feasible, convex optimization problem. The objective employed in this optimization problem yields a measure of distance (in a parameter set) to the power flow solution boundary. In practical terms, this distance is closely related to quantities that previous authors have proposed as voltage stability margins.  A typical margin is expressed in terms of the parameters of system loading (injected powers); here we additionally introduce a new margin in terms of the parameters of regulated bus voltages.


\end{abstract}


%
\IEEEpeerreviewmaketitle

\section{Introduction}
%
%
%
%


\IEEEPARstart{P}{ower} flow studies are the cornerstone of power system analysis and design. They are used in planning, operation, economic scheduling, transient stability, and contingency studies \cite{saadat}. The power flow equations used in these studies model the relationship between voltages and active and reactive power injections in a power system. The nonlinear power flow equations may not have any solutions (the power flow equations are said to be insolvable). That is, it is possible to choose a set of power injections for which no valid corresponding voltage profile exists. Practical cases that may fail to have a solution include long-range planning studies in which the studied system may not be able to support projected loads, and contingency studies for which the loss of one or more components may yield a network configuration that is similarly inoperable for the specified injections. This paper presents a practically computable sufficient condition, that when satisfied, rigorously classifies a specified case as unsolvable. This method also provides controlled voltage and power injection margins that characterize a distance to the power flow solvability boundary.



In engineering practice, large scale nonlinear power flow equations are typically solved using iterative numerical techniques, most commonly Newton-Raphson or its variants \cite{glover_sarma_overbye}, that rely on an initial guess of the solution voltage magnitudes and angles. An important limitation of these techniques is that they are only locally convergent. That is, they generally do not converge to a solution from an arbitrary initial guess \cite{saadat}, and may show very high sensitivity and highly complex behavior with respect to initial conditions for certain study cases.  It is well recognized that a power flow problem may in general have a very large number of solutions; for example, the work of \cite{baillieul_brynes} establishes cases for which the number of solutions grows faster than polynomial with respect to network size.  For cases having multiple solutions, each solution has a set of initial conditions that converges to that solution in Newton-Raphson iteration.  Characterization of Newton-Raphson regions of attraction was the subject of \cite{thorp1989}, which demonstrated cases for which the boundaries of these attractive sets were factual in nature.  So despite the fact that very large scale problems (10's or 100's of thousands of unknowns) are routinely solved in power engineering practice, it is important to recognize that as parameters move outside of routine operating ranges, the behavior of these equations can be highly complex, and failure of convergence for a Newton-Raphson-based commercial software package is far from a reliable indication that no solution exists.



The properties of the Newton-Raphson iteration guarantee (under suitable differentiability assumptions) that the iteration must converge to the solution for an initial condition selected in a sufficiently small neighborhood about that solution \cite{garcia_zangwill}.  However, when a selected initial condition (or some limited set of multiple initial conditions) fails to yield convergence, the user of a Newton-Raphson-based software package is left with an indeterminate outcome: does the specified problem have no solution, or has the initial condition(s) simply failed to fall within the attractive set of a solution that does exist?

Conditions to guarantee existence of solutions to the power flow equations has been an active topic of study. For example, \cite{DOE1} describes sufficient conditions for power flow solution existence. However, as sufficient conditions, these are often conservative: a solution may exist for a much larger range of operating points than satisfy the sufficient conditions. Other work on sufficient conditions for power flow solvability includes \cite{ilic1992}, which focuses on the decoupled (active power-voltage angle, reactive power-voltage magnitude) power flow model. Reference \cite{test_illcond} describes a modified Newton-Raphson iteration tailored to the type of ill-conditioning that can appear in power systems problems.  While convergence to a solution may be judged a constructive sufficient condition to demonstrate solvability, such approaches do not escape the fundamental limitations of a locally convergent iteration. In more recent work, \cite{necessary_stsl} provides two necessary conditions for saddle-node bifurcation based on lines reaching their static transfer stability limits; however, this work does not yet provide a test for power flow solvability or define a distance to the power flow solvability boundary.


A measure of the distance to the solvability boundary (the set of operating points where a solution exists, but small perturbations may result in the insolvability of the power flow equations \cite{hiskens2001}) is desirable to ensure that power systems are operated with security margins. If a solution does not exist for a specified set of power injections, a measure of the distance to the solvability boundary indicates how close the power flow equations are to having a solution. Existing work in this area uses a Newton-Raphson optimal multiplier approach \cite{optimal_multiplier} to find the voltage profile that yields the closest power injections to those specified \cite{overbye1994}, \cite{overbye1995}. Techniques for finding load margins comprised of the minimum distance to the power flow solvability boundary include \cite{dobson1993} and \cite{alvarado1994}. An algorithm that combines continuation and non-linear optimization techniques to either solve the power flow equations, when possible, or calculate a measure of power flow insolvability is presented in \cite{continuation_optimization_pf}. Reference \cite{interior_point_unsolvable} describes an optimization problem that applies interior point methods to minimize the load shedding necessary to obtain solvable power flow equations. The minimum amount of load shedding is used as a measure of power flow insolvability. Investigating the worst-case load shedding necessary for power flow solvability is also discussed in \cite{lesieutre2008} and \cite{feng1998}. Reference \cite{pf_insolvability_comparison} summarizes and compares some of these power flow insolvability measures.



In this paper, we present a sufficient condition under which a specified power flow problem is guaranteed to have no solution. The computation yields as a by-product controlled voltage and power injection margins to the power flow solvability boundary.


The sufficient condition for power flow insolvability is based on an optimization problem that includes a relaxation of certain equality constraints in the power flow equations. Specifically, in this optimization problem, the voltages at slack and PV buses are not fixed, but instead have a one-dimensional degree of freedom (i.e., they are allowed to change in constant proportion). In Section \ref{l:solution_existence_proof}, we prove that the extra degree of freedom guarantees the modified power flow equations to have at least one solution.  In an idealized lossless case, one may interpret this as follows: a sufficiently high voltage profile allows the system to meet any specified power injections. By continuity from the lossless case, we argue that this will continue to hold for modest losses, as is typical of models for bulk transmission. With the relaxed problem feasible for some (sufficiently high) voltage profile, we establish a non-empty feasible set for the optimization.


With a non-empty feasible set established, the optimization problem then seeks to minimize the slack bus voltage magnitude (using the one-degree-of-freedom in voltage profiles), subject to the active and reactive power injection constraints of the power flow equations. Importantly, we will show that a further relaxed version of this optimization problem is a convex semidefinite programming problem, and hence has a practically computable global minimum.  If the global minimum slack bus voltage obtained from this optimization problem is greater than the originally specified slack bus voltage, there can be no solution to the originally specified power flow equations. However, due to the nature of the relaxation, one may not draw a firm conclusion from the converse: if the minimum slack bus voltage is less than or equal to the specified slack bus voltage, the power flow equations may or may not be solvable. 

The ratio of the specified slack bus voltage to the minimum slack bus voltage gives a ``controlled voltage margin'' to the power flow solvability boundary. In a provably unsolvable case, this margin is the multiplicative factor by which the controlled voltages must be increased to allow the possibility of existence for a power flow solution. 

It is widely recognized that the power flow equations are quadratic in the complex voltage vector when these voltages are expressed in rectangular form.  Exploiting this fact, an analogous power injection margin can also be calculated; here the new, one degree of freedom introduced represents a constant power factor change in injections at each bus in proportion to the specified injections. When the power flow equations do not have a solution, the power injection margin provides the factor by which the power injections must be decreased to admit the possibility for power flow solution existence.

These margins are non-conservative upper bounds. Thus, for an insolvable set of specified values, a change in voltage by \emph{at least} the amount indicated by the voltage margin (or a change in power injections by at least the amount indicated by the power injection margin) is required for the power flow equations to be \emph{potentially} solvable. 




The dual of the optimization problem used in the sufficient condition can be written as a semidefinite program (SDP). The optimal power flow problem (i.e., finding the optimal operating point for a power system subject to physical and engineering constraints) was recently formulated as an SDP \cite{sdp_china, lavaei}. In prior work, the authors created an SDP formulation of the power flow equations in an attempt to calculate multiple solutions to these equations \cite{allerton2011}. In contrast to the non-convex primal optimization problem \cite{bernie_opfconvexity}, the feasible region of the dual problem formulated as an SDP is convex. The optimal objective value obtained from the dual SDP formulation is a lower bound on the objective function value used in the sufficient condition here. Thus, if the sufficient condition holds based on the lower bound from the dual SDP formulation, one can be assured the originally formulated power flow equations admit no solution.

The organization of this paper is as follows. In Section \ref{l:power_flow_equations_overview}, we give an overview of the power flow equations. In Section \ref{l:solution_existence_proof}, we provide the existence proof that shows the feasibility of the optimization problem used by the proposed condition. In Section \ref{l:cond1}, we describe the sufficient condition for power flow insolvability and define voltage and power injection margins. Numeric examples are provided in Section \ref{l:Examples1}. Section \ref{l:nonzero_gap} examines solutions with non-zero duality gap. We then conclude with a discussion of future work.


\section{Power Flow Equations Overview}
\label{l:power_flow_equations_overview}

The power flow equations describe the sinusoidal steady state equilibrium of a power network, and hence are formulated in terms of complex ``phasor'' representation of circuit quantities (see, for example, Ch. 9 of \cite{Chua_Desoer_Kuh}).  The underlying voltage-to-current relationships of the network are linear, but the nature of equipment in a power system is such that injected/demanded complex power at a bus (node) is typically specified, rather than current. The relation of interest is between the active and reactive power injected at each bus and the complex voltages at each bus, and hence the associated equations are nonlinear. Using the standard polar representation for complex voltages, and rectangular ``active/reactive'' representation of complex power, the power balance equations at bus $i$ are given by





\begin{align}
\label{P} P_{i} & = V_{i} \displaystyle\sum_{k=1}^n V_k \left( \mathbf{G}_{ik} \cos\left( \delta_i - \delta_k \right) + \mathbf{B}_{ik} \sin\left( \delta_i - \delta_k \right) \right) \\
\label{Q} Q_{i} & = V_{i} \displaystyle\sum_{k=1}^n V_k \left( \mathbf{G}_{ik} \sin\left( \delta_i - \delta_k \right) - \mathbf{B}_{ik} \cos\left( \delta_i - \delta_k \right) \right)
\end{align}

\noindent where $P_i$ and $Q_i$ are the active and reactive power injections, respectively, at bus $i$, $\mathbf{Y} = \mathbf{G} + j \mathbf{B}$ is the network admittance matrix, and $n$ is the number of buses in the system.



To represent typical behavior of equipment in the power system, each bus is classified as PQ, PV, or slack, according to the constraints imposed at that bus. PQ buses, which typically correspond to loads, treat $P_i$ and $Q_i$ as specified quantities, and enforce the active power \eqref{P} and reactive power \eqref{Q} equations at that bus. PV buses, which typically correspond to generators, specify a known voltage magnitude $V_i$ and active power injection $P_i$, and enforce only the active power equation \eqref{P}. The associated reactive power $Q_i$ may be computed as an ``output quantity,'' via \eqref{Q}. Finally, a single slack bus is selected, with its specified $V_{i}$ and $\delta_i$ (typically chosen to be $0^\circ$). The active power $P_i$ and reactive power $Q_i$ at the slack bus are determined from \eqref{P} and \eqref{Q}; network-wide conservation of complex power is thereby satisfied.

\section{Solution Existence Proof}
\label{l:solution_existence_proof}

The sufficient condition for power flow insolvability requires the evaluation of an optimization problem in which the feasible set is defined by a modified form of the power flow equations. The modification introduces one new degree of freedom, allowing voltage magnitudes at the slack and PV buses to vary; this variation is restricted to a one-degree-of-freedom ``ray,'' with all voltage magnitudes changing in constant proportion to their base-case values. We prove that the feasible space is non-empty for any lossless power system (i.e., a network model in which line conductances are all zero). Using standard results of basic circuit theory and continuity, we argue that the problem retains a non-empty feasible set when perturbed with small line conductances, as are typical in bulk transmission. 

The proof of solution existence may be outlined as follows. We first establish that a solution must exist for any lossless system with zero power injections. We then use the implicit function theorem to establish that solutions continue to exist for injections within small ball around zero. Hence, within this ball must exist a ray that aligns with the originally specified vector of non-zero power injections. We exploit the quadratic nature of the power flow equations to ``scale up'' voltage magnitudes along our one-degree-of-freedom, observing that the power injections must likewise move along the previously identified ray. It follows trivially that there will exist a scaling of voltages such that the specified power injections are realized, thereby yielding a solution to our modified form of the power flow equations.


\subsection{Existence of a Zero Power Injection Solution}
\label{l:zero_Pinj_exist}


Consider a generic lossless power system with all active and reactive power injections at PQ buses set to zero and all active power injections at PV buses set to zero. As our goal is accomplished if we can establish existence of one solution, we restrict attention to candidate solutions in which all buses have the same voltage angle of zero.

First, since zero power injection at a PQ bus implies zero nodal current injection, such buses have only branch admittances incident (i.e., from a circuit perspective, these are nodes with no independent source connected to them).  They can be eliminated from the network, and the network admittance matrix algebraically reduced via standard results of linear circuit theory. We generically assume that the reduced network does not result in any zero impedance lines.\footnote{Such a zero impedance line outcome can be eliminated by an arbitrarily small perturbation to the underlying line parameter data.}

Next, the substitution theorem \cite{Chua_Desoer_Kuh} guarantees that at any PV bus that has an associated non-zero reactive power injection, there must exist a shunt admittance of appropriate value such that, when substituted in place of the reactive injection, an identical solution for bus voltages is preserved.  The fact that the injections being replaced are purely reactive ensures that the associated admittances will be purely imaginary; i.e., susceptances only. 

With PQ buses eliminated and reactive injections at PV buses replaced by equivalent susceptances, the resulting network has the property that active and reactive power injections at all non-slack buses are identically zero. This allows us to write the remaining network constraints of interest as linear voltage/current relationships, as follows

\small
\begin{equation}\label{proof_VI_relationship}\left[\begin{array}{c} I_{\mathrm{slack}} \\ \hline \mathbf{0} \end{array}\right] = \left[\begin{array}{c|c}
j b_1 & j b_2 \\ \hline
j b_2^T & j \mathbf{B}_3 + j \mathrm{diag}\left(\Delta d\right)
\end{array}\right] \left[\begin{array}{c} 
V_{\mathrm{slack}} \\ \hline
V_{\mathrm{PV}}
\end{array}\right]
\end{equation}
\normalsize

\noindent where $\Delta d$ is is a vector of shunt element susceptances, $\mathrm{diag}\left(\Delta d \right)$ denotes the diagonal matrix with elements of $\Delta d$ on the diagonal, $\mathbf{B} = \left[ \begin{array}{c|c} b_1 & b_2 \\ \hline 
b_2^T & \mathbf{B}_3 \end{array} \right]$ is the bus susceptance matrix, and superscript $T$ indicates the transpose operator. $V_{slack}$ and $I_{slack}$ are the voltage and current injection at the slack bus, respectively, and $V_{PV}$ is the vector of PV bus voltages. Note that the lossless assumption implies that the network admittance matrix is purely imaginary.

Solving \eqref{proof_VI_relationship} for $\Delta d$ yields

\normalsize
\begin{equation} \label{DeltaD}
\Delta d = \left(\mathrm{diag}\left(V_{\mathrm{PV}} \right) \right)^{-1} \left(-b_2 V_{\mathrm{slack}} - \mathbf{B}_3 V_{\mathrm{PV}}\right)
\end{equation}
\normalsize

Because the voltage profile solution we seek is restricted to have the same voltage angle at all buses and a non-zero voltage magnitude at the slack bus, it follows that the voltage at every bus must be non-zero and $\mathrm{diag}\left(V_{\mathrm{PV}}\right)$ is invertible. Hence, for a lossless system under the assumptions specified, \eqref{DeltaD} yields a unique solution for the shunt susceptance values whose existence follows from the Substitution Theorem. 

Thus, the vector $\left[\begin{array}{c} 
V_{\mathrm{slack}} \\ \hline
V_{\mathrm{PV}}
\end{array}\right]$ provides a zero power injection solution to the reduced network that resulted from elimination of PQ buses; voltages at PQ buses can be trivially reconstructed. We conclude that any lossless system is guaranteed to have a zero power injection solution.

\begin{figure}[!th]
\centering
\includegraphics[totalheight=0.1\textheight]{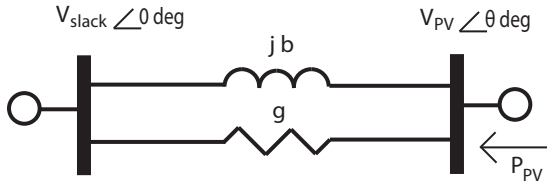}
\caption{Two-Bus System}
\label{twobussystemdiagram}
\end{figure}

To illustrate that this need not be the case for systems with large conductive elements in their bus admittance matrix (i.e., high transmission losses), consider the two-bus system with a slack bus and a PV bus shown in Fig. \ref{twobussystemdiagram}. 

The transmission line admittance is $g + j b$; note that in this admittance representation, the conductive term $g$ and the susceptance term $jb$ appear as parallel branch elements between the two buses. The voltage at the slack bus is denoted by $V_{slack}$, and the voltage at the PV bus is represented by $V_{PV}$ with angle $\theta$. The power injection at the PV bus is 

\begin{equation}\label{twobusequation} P_{PV} = g V_{PV}^2 - V_{PV} V_{slack} \left( g \cos\left(\theta\right) + b\sin\left(\theta\right) \right) \end{equation}

The two-bus system has a zero power injection solution for a given set of parameters $g$, $b$, $V_{PV}$, and $V_{slack}$ if a value of $\theta_0$ exists such that $P_{PV}\left( \theta_0 \right) = 0$. The existence of such a value of $\theta_0$ depends on the ratio of $V_{PV}$ to $V_{slack}$ and the ratio of $b$ to $g$. A zero power injection solution to this system exists when line resistances are small relative to line reactances and voltage magnitude differences are small. Specifically, the existence of a zero power injection solution for the system in Fig. \ref{twobussystemdiagram} requires

\begin{equation}\label{twobusreq} \left(\frac{V_{PV}}{V_{slack}}\right)^2 \leq 1 + \left(\frac{b}{g}\right)^2 \end{equation}

Since voltage magnitudes differences and line resistance to reactance ratios are typically small in realistic power systems, we expect that typical systems must have zero power injection solutions.  Consistent with this observation, all the IEEE power flow test cases of \cite{ieee_test_cases} have zero power injection solutions.  However, \eqref{twobusreq} confirms that the two-bus example will fail to have a zero injection solution when the conductance values relative to susceptance is sufficiently large; a zero injection solution fails to exist in the case of transmission parameters associated with very high losses. 


\subsection{Implicit Function Theorem}
\label{l:implicit_function_theorem}

We next apply the implicit function theorem \cite{hirsch_smale} at the zero power injection solution. Application of the implicit function theorem requires a non-singular Jacobian at the zero power injection solution. We therefore first investigate the Jacobian evaluated at this solution.

The Jacobian obtained using polar voltage coordinates at a zero power injection solution of a lossless system can be written as

\begin{equation}\label{lossless_jac}
\mathbf{J} = \begin{bmatrix}
\mathbf{J}_{11} & \mathbf{0} \\
\mathbf{0} & \mathbf{J}_{22}
\end{bmatrix}
\end{equation}

\noindent where

\begin{align}\label{lossless_jac11}
\mathbf{J}_{11} & = \frac{\partial P}{\partial \delta} = -\mathrm{diag}\left(V\right) \left( \mathbf{B}\, \mathrm{diag}\left( V \right) -  \mathrm{diag}\left(\mathbf{B} V \right)\right)\\
\label{lossless_jac22}
\mathbf{J}_{22} & = \frac{\partial Q}{\partial V} = -\mathrm{diag}\left(V \right) \mathbf{B} - \mathrm{diag}\left(\mathbf{B} V \right)
\end{align}

\noindent Since active and reactive power injections at the slack bus are unconstrained, the rows and columns corresponding to the slack bus are removed from both $\mathbf{J}_{11}$ and $\mathbf{J}_{22}$. Similarly, since the reactive power injections at PV buses are unconstrained, the rows and columns corresponding to PV buses are removed from $\mathbf{J}_{22}$. Note that both $\frac{\partial P}{\partial V}$ and $\frac{\partial Q}{\partial \delta}$ equal zero for the voltage profile with the same voltage angle at all buses corresponding to a zero power injection solution of a lossless system.


The implicit function theorem can be applied at a zero power injection solution so long as the Jacobian at this solution is non-singular. In the lossless case, this requires that $\mathbf{J}$ in \eqref{lossless_jac} is non-singular. We next show that the Jacobian for a lossless power system is non-singular at a zero power injection solution, provided that all lines are inductive and that the network is connected (i.e., no islands).

$\mathbf{J}$ in \eqref{lossless_jac} is non-singular if both $\mathbf{J}_{11}$ and $\mathbf{J}_{22}$ are non-singular. The matrix $\mathrm{diag}\left(\mathbf{B} V \right)$ in $\mathbf{J}_{22}$ is equivalent to $\mathrm{diag}\left( I \right)$, where $I$ is the vector of current injections. Since all rows and columns in $\mathbf{J}_{22}$ correspond to PQ buses with zero current injections, this term is zero. The matrix $\mathrm{diag}\left(V\right)$ is non-singular since all voltages are non-zero for the voltage profile with the same voltage angle at all buses. With the slack bus row and column removed, $\mathbf{B}$ is non-singular for a connected power system with inductive lines. Thus, $\mathbf{J}_{22}$ is non-singular. 

Since $\mathrm{diag}\left(V\right)$ is non-singular, $\mathbf{J}_{11}$ is non-singular if 


\small
\begin{align}\label{J11_nonsingular}\nonumber & \left( \mathbf{B}\,\mathrm{diag}\left( V \right) -  \mathrm{diag}\left(\mathbf{B} V \right)\right) =  \\
& \begin{bmatrix}
\mathbf{B}_{12}V_2 + \ldots	+ \mathbf{B}_{1n}V_n & \cdots & -\mathbf{B}_{1n}V_n \\
\vdots & \ddots & \vdots \\
-\mathbf{B}_{n1}V_1  & \cdots & \mathbf{B}_{n1}V_1 + \ldots + \mathbf{B}_{n\left(n-1\right)}V_{\left(n-1\right)} \\
\end{bmatrix}
\end{align}
\normalsize

\noindent is non-singular. Note that the diagonal elements in \eqref{J11_nonsingular} are the negative of the sum of the off-diagonal elements in the corresponding row. Under the assumption of inductive lines such that all off-diagonal elements are negative, this matrix has weak diagonal dominance. With the slack bus row and column removed, the remaining matrix has at least one row where the diagonal element is strictly greater than the sum of the off-diagonal elements (i.e., strict diagonal dominance exists for this row). Since the power system associated with the $\mathbf{B}$ matrix is connected (i.e., no islands), the digraph associated with the matrix in \eqref{J11_nonsingular} is strongly connected. This implies that the matrix is irreducible \cite{horn1985}. Since the matrix is irreducible, weakly diagonally dominant, and has at least one row with strict diagonal dominance, the matrix is irreducibly diagonally dominant. By the Levy--Desplanques theorem, the matrix is non-singular \cite{horn1985}. Thus, $\mathbf{J}_{11}$ is non-singular. This proves that the Jacobian for a connected, lossless system at a zero power injection solution is invertible, under the assumption of inductive lines.


Although the assumptions of lossless systems and inductive lines are required for the above proof, non-singularity of the Jacobian at a zero power injection solution generically holds for more general systems (e.g., lossless systems with some capacitive lines and lossy systems). A singular Jacobian would imply marginal stability at the zero power injection solution with multiple solutions coalescing at a bifurcation point. There is no reason to expect this to occur at a zero power injection solution. Zero power injection solutions for all IEEE power flow test cases \cite{ieee_test_cases} have non-singular Jacobians at zero power injection solutions.

If the Jacobian of the power flow equations is non-singular at the zero power injection solution, the implicit function theorem indicates that a solution must persist for all power injections in a small ball around the zero power injection.  Thus, there exists some voltage magnitude and angle perturbation $\Delta V\angle \Delta \delta$ such that 

\begin{equation}\label{impfunct} f\left( V + \Delta V \angle \Delta \delta \right) = \Delta P + j \Delta Q \end{equation}

\noindent for any small $\Delta P$, $\Delta Q$, where $V$ is the voltage profile for the zero power injection solution, $\Delta P$ and $\Delta Q$ are small perturbations to the active and reactive power injections, and $f$ represents the power flow equations relating the voltages and power injections.

\subsection{Scaling Up Voltages}
\label{l:scaling_up_voltages}

We complete the solution existence proof by expanding the small ball around the zero power injection solution to obtain a voltage profile that yields the originally specified power injections. Since the power flow equations are quadratic in voltage magnitudes $V$, scaling all voltage magnitudes also scales the power injections. That is, scaling the voltage magnitudes in \eqref{impfunct} by the scalar $\beta$ gives

\begin{equation}\label{scaleV} f\left(\beta \left( V + \Delta V \angle \Delta \delta \right) \right) = \beta^2 \left(\Delta P + j \Delta Q\right) \end{equation}

Choose a $\Delta P + j \Delta Q$ that is in the direction of the specified power injections and obtain a corresponding voltage profile $V + \Delta V\angle \Delta \delta$. Then increase $\beta$ until the power injections given by $f\left(\beta \left( V + \Delta V\angle \Delta \delta\right)\right)$ match the specified power injections. The voltage profile $\beta \left( V + \Delta V\angle \Delta \delta\right)$ then yields the specified power injections. 



\section{Sufficient Condition for Power Flow Insolvability}
\label{l:cond1}

\subsection{Condition Description}
\label{l:cond1description}

The proof in Section \ref{l:solution_existence_proof} shows that there exists a voltage profile satisfying the power injection equations. We develop the sufficient condition for power flow insolvability by determining whether any such voltage profile could match the specified slack bus and PV bus voltages. No solution to the power flow equations exists if it is impossible to obtain a voltage profile that yields the specified power injections while also matching the specified voltage magnitudes at slack and PV buses. 

One way to determine if a valid voltage profile exists is to find the voltage profile with the lowest possible slack bus voltage. If the minimum possible slack bus voltage is greater than the specified slack bus voltage, no voltage profile will satisfy the power flow equations and thus the power flow equations are insolvable. If, conversely, the minimum slack bus voltage is less than the specified slack bus voltage, a power flow solution may exist.

This condition thus indicates that no power flow solution exists when the minimum slack bus voltage obtainable while satisfying the power injection equations (with PV bus voltage magnitudes scaled proportionally) is greater than the specified slack bus voltage magnitude. An optimization problem with objective function minimizing the slack bus voltage and constraints on power injections and PV bus voltage magnitudes, as shown in \eqref{cond1}, is used to evaluate this condition.


\begin{subequations}
\small
\label{cond1}
\begin{align}\label{cond1obj}
& \mathrm{min}\quad V_{\mathrm{slack}} \\\nonumber
& \mathrm{subject\; to} \\\nonumber
& P_{k} = V_{k} \displaystyle\sum_{i=1}^n V_i \left( \mathbf{G}_{ik} \cos\left( \delta_k - \delta_i \right) + \mathbf{B}_{ik} \sin\left( \delta_k - \delta_i \right) \right) \\[-8pt] \label{cond1P}
& \qquad\qquad\qquad\qquad\qquad\qquad\qquad\quad \forall\, k \in \left\lbrace\mathcal{PQ,\, PV} \right\rbrace \\\nonumber
& Q_{k} = V_{k} \displaystyle\sum_{i=1}^n V_i \left( \mathbf{G}_{ik} \sin\left( \delta_k - \delta_i \right) - \mathbf{B}_{ik} \cos\left( \delta_k - \delta_i \right) \right) \\[-8pt] \label{cond1Q}
& \qquad\qquad\qquad\qquad\qquad\qquad\qquad\quad \forall\, k \in \mathcal{PQ} \\[5pt] \label{cond1V}
& V_k  = \alpha_k V_{\mathrm{slack}} \qquad\qquad\qquad\qquad\quad\,	 \forall\, k \in \mathcal{PV} 
\end{align}
\end{subequations}

\noindent where $\mathcal{PQ}$ is the set of PQ buses, $\mathcal{PV}$ is the set of PV buses, and $V_{slack}$ is the slack bus voltage magnitude. $\alpha_k$ represents the specified ratio of the PV bus $k$ and slack bus voltage magnitudes. The minimum achievable slack bus voltage (i.e., the optimal objective value of \eqref{cond1}) is denoted as $V_{slack}^{min}$.

Note that \eqref{cond1} is feasible for lossless systems and expected to be feasible for most practical power systems as argued in Section \ref{l:solution_existence_proof}. The power flow insolvability condition can therefore be evaluated.

The optimization problem \eqref{cond1} is in general non-convex \cite{bernie_opfconvexity}, and hence solution for a global optimum is not assured. A global minimum is required in order to ensure the validity of the sufficient condition on power flow solution non-existence. We therefore formulate in \eqref{cond1dual} the semidefinite dual of \eqref{cond1}. SDP algorithms can assure that we find the global solution to the convex dual formulation \eqref{cond1dual}.

\begin{subequations}
\small
\label{cond1dual}
\begin{align}  & \mathrm{max}  \displaystyle\sum_{k \in \left\lbrace \mathcal{PQ},\, \mathcal{PV} \right\rbrace } \!\!\!\!\!\!\! \left( \lambda_k P_{k} \right) + \displaystyle\sum_{k \in \mathcal{PQ} }\left( \gamma_k Q_{k} \right) \\ \nonumber
& \mathrm{subject\; to} \\\nonumber
& \mathbf{A}\left(\lambda, \gamma, \mu \right) = \left[ \mathbf{M}_{\mathrm{slack}} - \displaystyle\sum_{k \in \mathcal{PQ}} \left( \lambda_k \mathbf{Y}_k + \gamma_k \mathbf{\bar{Y}}_k \right) \right. \\ \label{cond1dualA}
& \qquad\qquad \left. - \displaystyle\sum_{k \in \mathcal{PV}} \left( \lambda_k \mathbf{Y}_k + \mu_k \left( \mathbf{M}_k - \alpha_k \mathbf{M}_{\mathrm{slack}} \right) \right) \right] \succeq 0 \end{align}
\end{subequations}
\normalsize

\noindent where $\lambda_k$, $\gamma_k$, and $\mu_k$ are the Lagrange multipliers for active power (equation \eqref{cond1P}), reactive power (equation \eqref{cond1Q}), and PV bus voltage magnitude ratio (equation \eqref{cond1V}) constraints, respectively, associated with bus $k$. The symbol $\succeq$ indicates that the corresponding matrix is constrained to be positive semidefinite. The maximum lower bound on the minimum achievable slack bus voltage (i.e., the optimal objective value of \eqref{cond1dual}) is denoted as $\underline{V}_{slack}^{min}$.

Matrices employed in \eqref{cond1dual} are defined as

\small
\begin{align}\label{Yk}\mathbf{Y}_k & = \frac{1}{2} \begin{bmatrix}
\operatorname{Re} \left(Y_k + Y_k^T\right) & \operatorname{Im} \left(Y_k^T - Y_k\right) \\
\operatorname{Im} \left(Y_k - Y_k^T\right) & \operatorname{Re} \left(Y_k\ + Y_k^T\right)
\end{bmatrix} \\
\label{Ykbar}\mathbf{\bar{Y}}_k & = -\frac{1}{2} \begin{bmatrix}
\operatorname{Im} \left(Y_k + Y_k^T\right) & \operatorname{Re} \left(Y_k - Y_k^T\right) \\
\operatorname{Re} \left(Y_k^T - Y_k\right) & \operatorname{Im} \left(Y_k + Y_k^T\right)
\end{bmatrix} \\
\label{Mk}\mathbf{M}_k & = \begin{bmatrix}
e_k e_k^T & \mathbf{0} \\
\mathbf{0} & e_k e_k^T
\end{bmatrix} \end{align}
\normalsize

\noindent where $e_k$ denotes the $k^{th}$ standard basis vector in $\mathbb{R}^n$ and the matrix $Y_k = e_k e_k^T \mathbf{Y}$. Notation is adopted from \cite{lavaei}. 

Note that the dual formulation \eqref{cond1dual} is always feasible since the point $\lambda_i = 0$, $\gamma_i = 0$, $\mu_i = 0$  for all $i$ implies $\mathbf{A} = \mathbf{M}_{slack} \succeq 0$.

The semidefinite dual formulation \eqref{cond1dual} provides a lower bound on the minimum slack bus voltage in \eqref{cond1}. No solution to the power flow equations exists if the lower bound from \eqref{cond1dual} is greater than the specified slack bus voltage. That is, the condition 

\begin{equation}\label{cond1g}\underline{V}_{slack}^{min} > V_{0}\end{equation}

\noindent where $V_0$ is the specified slack bus voltage, is a sufficient condition for insolvability of the power flow equations. Note that this formulation does not use any rank relaxations or enforce any requirements on matrix ranks; the solution to the convex problem \eqref{cond1dual} is only used as a lower bound on \eqref{cond1}.

The converse condition does not necessarily hold: the power flow equations may not have a solution even if 

\begin{equation}\label{cond1leq}\underline{V}_{slack}^{min} \leq V_{0}\end{equation}

\noindent Thus, \eqref{cond1leq} is a necessary, but not sufficient, condition for power flow solvability. However, satisfaction of \eqref{cond1leq} is expected to often predict the existence of a power flow solution. 

If the $\mathbf{A}$ matrix in \eqref{cond1dualA} has a nullspace with rank less than or equal to two, a solution of the power flow equations with slack bus voltage equal to $\underline{V}_{slack}^{min}$ (and PV bus voltage magnitudes scaled proportionally) can be obtained (see \cite{lavaei} for further details). If a solution to the power flow equations with slack bus voltage equal to $V_0$ does not exist, the solution with lower slack bus voltage must disappear as the controlled voltages increase. The disappearance of a solution due to increasing controlled voltages does not typically occur. Thus, satisfaction of \eqref{cond1leq} by a solution to \eqref{cond1dual} with $\mathrm{rank}\left(\mathrm{null}\left(\mathbf{A}\right)\right) \leq 2$ is a strong indicator of power flow solution existence.

\subsection{Controlled Voltage Margin}
\label{l:voltage_margin}

The sufficient condition \eqref{cond1g} is binary: the specified power flow equations either cannot have a solution or may have a solution. We next show how the sufficient condition can be interpreted to give a measure of the \emph{degree} of solvability. We develop a measure of the distance to the power flow solvability boundary, which we define as the set of solvable power injections where all solutions may vanish under small perturbations. Since operating a power system far from the power flow solvability boundary is required to ensure stability, a measure of the distance to the solvability boundary is useful. A measure of the distance to the solvability boundary also indicates how close insolvable power flow equations are to solvability.

We introduce a controlled voltage margin measure $\sigma$ for the distance to the power flow solvability boundary. The controlled voltage margin is defined as the ratio between the specified slack bus voltage and the lower bound on the minimum slack bus voltage obtained from \eqref{cond1dual}.

\begin{equation}\label{Vmargin} \sigma = \frac{V_{0}}{\underline{V}_{slack}^{min}} \end{equation}

$\sigma$ is an upper (non-conservative) bound of the distance to the power flow solvability boundary. For a solvable set of power flow equations, we are guaranteed to be at or beyond the solvability boundary if the specified slack bus voltage decreases by the factor $\sigma$. For insolvable power flow equations, increasing the slack bus voltage magnitude (with proportional increases in PV bus voltage magnitudes) by at least a factor of $\frac{1}{\sigma}$ (without changing the power injections) is required for solvability.

The sufficient condition can be rewritten in terms of the voltage margin: $\sigma < 1$ is a sufficient condition for power flow insolvability.


\subsection{Power Injection Margin}
\label{l:Pinj_margin}

Another measure of the distance to the power flow solvability boundary is given in terms of power injections. The power injection margin is a measure of how large of a change in the power injections in a certain profile is required for the power injections to be on the solvability boundary. We consider the profile where power injections are uniformly changed at each bus in order to take advantage of the quadratic nature of the optimization problem \eqref{cond1} in the sufficient condition. (The impact of non-uniform changes in power injections is future work discussed in Section \ref{l:conclusion}.) The quadratic property that we exploit can be written as

\begin{equation}\label{pinj_margin} h \left( \eta \left( P + j Q \right) \right) = \eta \left( \underline{V}_{slack}^{min} \right)^2 \end{equation}

\noindent where $P$ and $Q$ are vectors of the active and reactive power injection at each bus, $h$ is the function representing optimization problem \eqref{cond1} relating the minimum slack bus voltage to the power injections, and $\eta$ is a scalar.

\eqref{pinj_margin} describes the linear relationship between the square of the voltages and the power injections. This relationship is evident from \eqref{cond1P} and \eqref{cond1Q}: scaling all voltages by $\sqrt{\eta}$ scales the active and reactive power injections by $\eta$.

To develop the power injection margin, uniformly scale the power injections until the sufficient condition \eqref{cond1g} indicates that the power injections are (at least) on the solvability boundary.

\begin{equation}\label{pinj_margin_def}\eta \left( \underline{V}_{slack}^{min} \right)^2  = \left(V_{0}\right)^2 \end{equation}

The power injection margin $\eta$ corresponding to the condition in \eqref{pinj_margin_def} gives an upper, non-conservative bound of the distance to the solvability boundary in the direction of uniformly increasing power injections. For a solvable set of power injections, the largest proportional increase in power injections at each bus while potentially maintaining solvability is a factor of $\eta$. For an insolvable set of power injections, a proportional change of all power injections by at least $\eta$ is required for a solution to be possible.

Note that the power injection margin can be rewritten in terms of the voltage margin.

\begin{equation}\label{eta_sigma} \eta = \left(\sigma \right)^2 \end{equation}

The sufficient condition for power flow insolvability can be rewritten in terms of the power injection margin: $\eta < 1$ is a sufficient condition for power flow insolvability.

The power injection margin can alternatively be calculated using a different optimization problem that directly maximizes power injections at constant power factor with fixed slack and PV bus voltage magnitudes. This gives the same power injection margins as \eqref{pinj_margin_def} in all test cases. This approach lacks a feasibility proof and does not have any advantages over \eqref{pinj_margin_def}.

\section{Numeric Example}
\label{l:Examples1}

We next apply the sufficient condition for power flow solvability to the IEEE 14-bus and IEEE 118-bus systems \cite{ieee_test_cases} using optimization codes YALMIP \cite{yalmip} and SeDuMi \cite{sedumi}. The power injections are uniformly increased at each bus at constant power factor until the sufficient condition indicates that no solutions are possible. The sufficient condition results are compared to power flow solution attempts by a Newton-Raphson algorithm.

\subsection{IEEE 14-Bus System Results}
\label{l:14bus_cond1}

Results from applying the sufficient condition to the IEEE 14-bus system are given in Table \ref{14busresults}. The specified slack bus voltage is $V_0 = 1.0600$ per unit.

\begin{table}[ht]
\centering
\begin{tabular}{|c|c|c|}
\hline 
\textbf{Injection Multiplier} & \textbf{NR Converged} & $\underline{V}_{slack}^{min}$  \\ \hline\hline	
1.000 & Yes  & 0.5261\\ \hline
2.000 & Yes  & 0.7440\\ \hline
3.000 & Yes  & 0.9112\\ \hline
4.000 & Yes  & 1.0522\\ \hline
4.010 & Yes  & 1.0535\\ \hline
4.020 & Yes  & 1.0548\\ \hline
4.030 & Yes  & 1.0561\\ \hline
4.040 & Yes  & 1.0575\\ \hline
4.050 & Yes  & 1.0588\\ \hline
4.055 & Yes  & 1.0594\\ \hline
4.056 & Yes  & 1.0595\\ \hline
4.057 & Yes  & 1.0597\\ \hline
4.058 & Yes  & 1.0598\\ \hline
\textbf{4.059} & \textbf{Yes} & \textbf{1.0599}\\ \hline\hline
\textbf{4.060} & \textbf{No}  & \textbf{1.0601}\\ \hline
4.061 & No  & 1.0602\\ \hline
4.062 & No  & 1.0603\\ \hline
4.063 & No  & 1.0605\\ \hline
4.064 & No  & 1.0606\\ \hline
4.065 & No  & 1.0607\\ \hline
5.000 & No  & 1.1764\\ \hline
\end{tabular}
\caption{Solvability Condition Results For IEEE 14-Bus System}
\label{14busresults}
\end{table}
\normalsize

The originally specified active and reactive power injections are increased uniformly at each bus. The first column of Table \ref{14busresults} lists the multiple by which the injections are increased. No power flow solutions exist after a sufficiently large increase (approximately 4.060 for this example). Note that the injection multiplier given in the first column does not change at a constant rate but rather focuses on the region near power flow solution non-existence. 

The second column indicates whether a Newton-Raphson solver converged to a solution at the corresponding loading. In order to increase the likelihood of convergence, the Newton-Raphson solver was initialized at each injection multiplier with the solution from the previous injection multiplier and a large number of Newton-Raphson iterations were allowed.

The third column provides the lower bound on the minimum slack bus voltage in per unit obtained from \eqref{cond1dual}. In order to evaluate the sufficient condition for power flow insolvability at each injection multiplier, the value in this column is compared to the specified slack bus voltage of 1.06 per unit. If the value in the third column is greater than 1.06, the sufficient condition indicates that no power flow solutions exist. These results show agreement between Newton-Raphson convergence and the sufficient condition; a power flow solution was found for all injection multipliers where the sufficient condition indicated that a solution was possible (observe that both $\underline{V}_{slack}^{min}$ is just greater than 1.06 and no solution is found by the Newton-Raphson solver at an injection multiplier of 4.060). 

The existence of a solution for all power injections that satisfy \eqref{cond1leq} is expected since the $\mathbf{A}$ matrix in \eqref{cond1dualA} has a nullspace with rank two. This need not always be the case. In \ref{l:118bus_cond1}, we investigate an example with $\mathrm{rank}\left(\mathrm{null}\left(\mathbf{A}\right)\right) = 4$ where no solution was found for some power injections even though the condition \eqref{cond1leq} indicated that a solution was possible.

\begin{figure}[t]	
\centering
\includegraphics[totalheight=0.21\textheight]{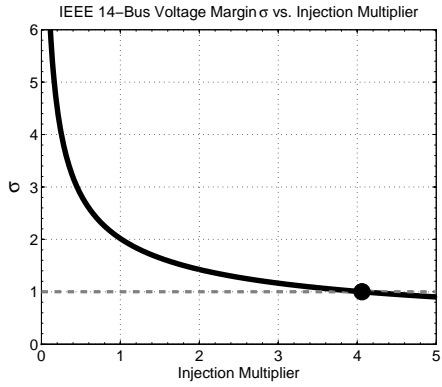}
\caption{IEEE 14-bus Voltage Margin}
\label{f:vmargin_14}
\end{figure}

We next use the IEEE 14-bus system example to demonstrate the voltage and power injection margins. In Fig. \ref{f:vmargin_14}, the voltage margin $\sigma$ is plotted versus the injection multiplier. The voltage margin decreases as power injections increase. The voltage margin crosses one at an injection multiplier of 4.0595, indicating that no power flow solution can exist for larger power injections. Beyond this point, the voltage margin provides the minimum increase in the slack bus voltage (with corresponding proportional voltage increases at all PV buses) required in order for a power flow solution to possibly exist. 

\begin{figure}[tb]
\centering
\subfloat[Nominal and Low Slack Bus Voltage PV Curves]{\includegraphics[totalheight=0.21\textheight]{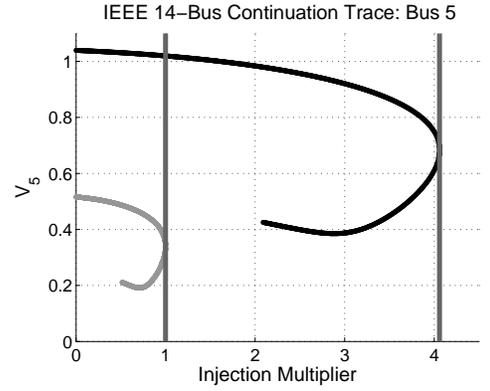}\label{f:pv_lv_14}} \hfill
\subfloat[Nominal and High Slack Bus Voltage PV Curves]{\includegraphics[totalheight=0.21\textheight]{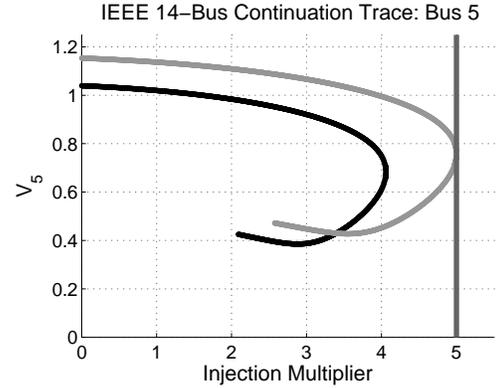}\label{f:pv_hv_14}} \\
\caption{IEEE 14-Bus System PV Curves}
\label{f:pv_14}
\end{figure}

In Fig. \ref{f:pv_14}, we examine the power versus voltage (PV) curves for the high-voltage, stable solution to the IEEE 14-bus system. These curves, which were plotted using continuation techniques \cite{cpf}, show how the solution voltages change with proportional increases in power injections at all buses. The plots show the voltage at the arbitrarily selected PQ bus five. (Plotting the voltage at a PQ bus is required since voltage magnitudes at slack and PV buses are fixed.) The PV curve using the nominal slack and PV bus voltages is shown in black.


Evaluating the optimization problem \eqref{cond1dual} at an injection multiplier of one gives a $\underline{V}_{slack}^{min} = 0.5261$. The voltage margin is $\sigma = \frac{1.0600}{0.5261} = 2.0148$ per unit. Thus, no solution can exist if the slack bus voltage is reduced by more than a factor of 2.0148 (with all PV bus voltages reduced proportionally). The grey PV curve in Fig. \ref{f:pv_lv_14} is obtained when the voltages are thus reduced. This curve shows that with these reduced voltages, there is the single solution is on the power flow solvability boundary; no solutions exist after any further increase in the injection multiplier. Thus, the voltage margin accurately indicates the distance to power flow insolvability.

The solution to the optimization problem \eqref{cond1dual} also enables determination of the power injection margin $\eta$. Solving \eqref{pinj_margin_def} yields $\eta = \left(\frac{1.0600}{0.5261}\right)^2 = 4.0595$. Thus, the power injections can be increased uniformly by a factor of 4.0595 until the sufficient condition indicates that no power flow solutions are possible. The black PV curve associated with the nominal voltages in Fig. \ref{f:pv_lv_14} corroborates this assertion: a power flow solution exists for all power injection multipliers less than 4.0595, but no solution exists beyond this power injection multiplier.

The voltage and power injection margins can also be used to investigate insolvable power injections. Assume that we would like to consider operation at a power injection multiplier equal to five. Evaluating the optimization problem \eqref{cond1dual} at a power injection multiplier of five gives $\underline{V}_{slack}^{min} = 1.1764$. Note that \eqref{pinj_margin_def} implies that knowledge of $\underline{V}_{slack}^{min}$ at a power injection multiplier of one allows the direct calculation $\underline{V}_{slack}^{min}$ at a power injection multiplier of five: 

\small
\begin{align} \nonumber \left. \underline{V}_{slack}^{min} \right|_{\mathrm{Inj\, Mult = 5}} & = \sqrt{\eta} \left. \underline{V}_{slack}^{min} \right|_{\mathrm{Inj\, Mult = 1}}\\
& = 0.5261 \cdot \sqrt{5}  = 1.1764\mathrm{\, per\, unit}\end{align}
\normalsize

The voltage margin at a power injection multiplier of five is $\sigma = \frac{1.06}{1.1764} = 0.9011$. $\sigma < 1$ indicates that there is no solution at a power injection multiplier of five. To potentially achieve a power flow solution, the slack bus voltage must increase by at least a factor of $\frac{1}{0.9011} = 1.1098$ (with corresponding proportional increases in all PV bus voltages). The grey PV curve in Fig. \ref{f:pv_hv_14} has the voltages thus increased. Observe that increasing the voltages in this way allows for a power flow solution on the power flow solvability boundary for a power injection multiplier of five.

The power injection margin $\eta$ can also be calculated at a power injection multiplier of five using \eqref{pinj_margin_def}.

\small 
\begin{align}\nonumber \eta & =  \left( \frac{V_{0}}{\left. \underline{V}_{slack}^{min} \right|_{\mathrm{Inj\, Mult = 5}}} \right)^2 \\
& = \left(\frac{1.0600}{1.1764} \right)^2 = 0.8119 \end{align}
\normalsize

$\eta < 1$ implies that no solution exists at a power injection multiplier of five. The power injection margin also indicates that no solution can exist for power injection multipliers greater than $ 0.8119 \cdot 5 = 4.0595$. This corresponds to the ``nose'' point of the black (nominal) PV curve in Fig. \ref{f:pv_hv_14}.

\subsection{IEEE 118-Bus System Results}
\label{l:118bus_cond1}

Results from applying the sufficient condition to the IEEE 118-bus system are given in Table \ref{118busresults}. The data are arranged in the same manner as in Table \ref{14busresults}. The specified slack bus voltage is $V_0 = 1.0350$ per unit. 


The results can be categorized intro three regions: small power injections where the sufficient condition indicates that a solution is possible and a solution is indeed found using a Newton-Raphson solver, larger power injections where the sufficient condition indicates that a solution is possible but no solutions are found, and yet larger power injections where the sufficient condition indicates that no solutions are possible and no solutions are found.

The rank of the nullspace of the $\mathbf{A}$ matrix in \eqref{cond1dualA} for the IEEE 118-bus system was four. Therefore, the expectation that satisfaction of \eqref{cond1leq} will result in power flow solvability may not hold. Correspondingly, these results emphasize the fact that \eqref{cond1g} is a \emph{sufficient} condition for power flow insolvability. Specifically, a power flow solution was not found for injection multipliers greater than 3.18, even though $\underline{V}_{slack}^{min}$ is less than the specified slack bus voltage until an injection multiplier of 3.27. A continuation power flow indicates that the high-voltage solution bifurcates at a power injection multiplier of 3.185, so it is likely that no solutions exist after this point. No solutions are found at injection multipliers larger than 3.27 where the sufficient condition indicates that no solutions are possible.

\begin{table}[t]
\centering
\begin{tabular}{|c|c|c|}
\hline 
\textbf{Injection Multiplier} & \textbf{NR Converged} & $\underline{V}_{slack}^{min}$  \\ \hline\hline	
1.00 & Yes & 0.5724\\ \hline
1.50 & Yes & 0.7010\\ \hline
2.00 & Yes & 0.8095\\ \hline
2.50 & Yes & 0.9050\\ \hline
3.00 & Yes & 0.9914\\ \hline
3.15 & Yes & 1.0159\\ \hline
3.16 & Yes & 1.0175\\ \hline
3.17 & Yes & 1.0191\\ \hline
3.18 & Yes & 1.0207\\ \hline\hline
3.19 & No & 1.0223\\ \hline
3.20 & No & 1.0239\\ \hline
3.21 & No & 1.0255\\ \hline
3.22 & No & 1.0271\\ \hline
3.23 & No & 1.0287\\ \hline
3.24 & No & 1.0303\\ \hline
3.25 & No & 1.0319\\ \hline
3.26 & No & 1.0335\\ \hline\hline
3.27 & No & 1.0351\\ \hline
3.28 & No & 1.0366\\ \hline
3.29 & No & 1.0382\\ \hline
4.00 & No & 1.1448\\ \hline
\end{tabular}
\caption{Insolvability Condition Results for IEEE 118-Bus System}
\label{118busresults}
\end{table}
\normalsize

\begin{figure}[t]
\centering
\includegraphics[totalheight=0.25\textheight]{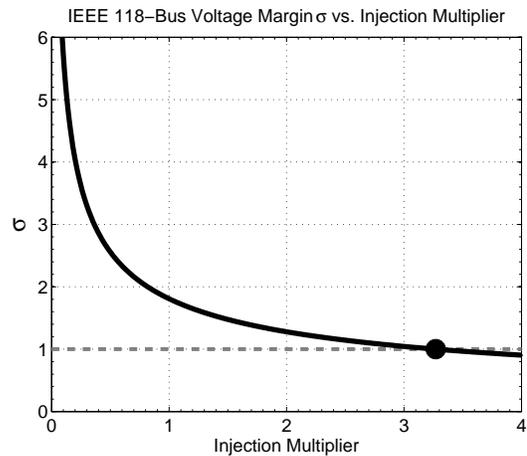}
\caption{IEEE 118-bus Voltage Margin}
\label{f:vmargin_118}
\end{figure}

We next use the IEEE 118-bus system example to demonstrate the voltage and power injection margins. In Fig. \ref{f:vmargin_118}, the voltage margin $\sigma$ is plotted versus the injection multiplier. Similar to Fig. \ref{f:vmargin_14}, the voltage margin decreases as power injections increase. The voltage margin crosses one at an injection multiplier of $3.2695$, indicating that no power flow solution can exist for larger power injections. For larger power injections, the voltage margin provides the minimum increase in the slack bus voltage (with corresponding proportional voltage increases at all PV buses) that is required in order for a power flow solution to possibly exist.

\begin{figure}[!t]
\subfloat[Nominal and Low Slack Bus Voltage PV Curves]{\includegraphics[totalheight=0.25\textheight]{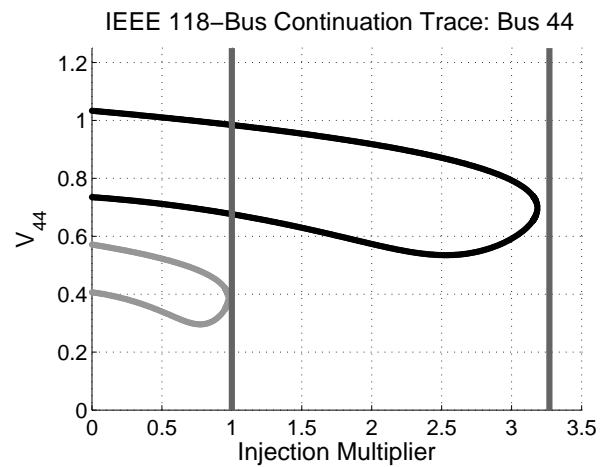}\label{f:pv_lv_118}} \hfill
\subfloat[Nominal and High Slack Bus Voltage PV Curves]{\includegraphics[totalheight=0.25\textheight]{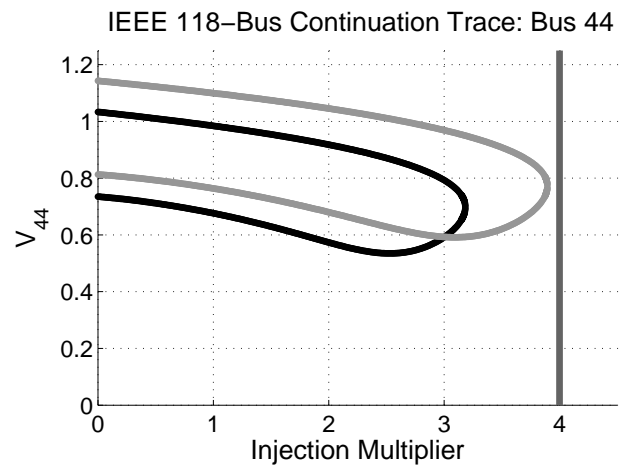}\label{f:pv_hv_118}} \\
\caption{IEEE 118-Bus System PV Curves}
\label{f:pv_118}
\end{figure}

In Fig. \ref{f:pv_118}, we examine the PV curves for the high-voltage, stable solution to the IEEE 118-bus system. The plots show the voltage at the arbitrarily selected PQ bus 44 (plotting the voltage at a PQ bus is required since the voltages at slack and PV buses are fixed). The PV curve using the nominal slack and PV bus voltages is shown in black.


Evaluating the optimization problem \eqref{cond1dual} associated with the sufficient condition at an injection multiplier of one gives a $\underline{V}_{slack}^{min} = 0.5724$. The voltage margin is $\sigma = \frac{1.0350}{0.5724} = 1.8082$. Thus, no solution can exist if the slack bus voltage is reduced by more than a factor of 1.8082 (with all PV bus voltages reduced proportionally). The grey PV curve in Fig. \ref{f:pv_lv_118} is obtained when the voltages are thus reduced. Although no solutions exist for the grey PV curve at injection multipliers larger than one, there are also injection multipliers slightly less than one for which no solutions are found with continuation techniques. This reinforces the fact that the voltage margin is an \emph{upper bound} on the distance to the solvability boundary.

The solution to the optimization problem \eqref{cond1dual} also enables determination of the power injection margin $\eta$. Solving \eqref{pinj_margin_def} yields $\eta = \left(\frac{1.0350}{0.5724}\right)^2 = 3.2695$. Thus, the power injections can be increased uniformly by a factor of 3.2695 until the sufficient condition indicates that no power flow solutions are possible. This is also an upper bound on the distance to the solvability boundary: as indicated by the sufficient condition, the black PV curve associated with the nominal voltages in Fig. \ref{f:pv_lv_118} has no solutions for power injection multipliers larger than 3.2695, but also appears to have no solutions for some values of power injection multipliers below 3.2695. (It is possible, but unlikely, that a PV curve associated with a different solution may exist at injection multipliers between the ``nose'' of the PV curve associated with the high-voltage solution at 3.1840 and the value of 3.2695 from the sufficient condition.)

The voltage and power injection margins can also be used to investigate insolvable power injections. Assume that we would like to consider operation at a power injection multiplier equal to four. Evaluating the optimization problem \eqref{cond1dual} at a power injection multiplier of four gives $\underline{V}_{slack}^{min} = 1.1448$. Note that \eqref{pinj_margin_def} implies that knowledge of $\underline{V}_{slack}^{min}$ at a power injection multiplier of one allows the direct calculation $\underline{V}_{slack}^{min}$ at a power injection multiplier of four: 

\begin{align}\nonumber \left. \underline{V}_{slack}^{min} \right|_{\mathrm{Inj\, Mult = 4}} & = \sqrt{\eta} \left. \underline{V}_{slack}^{min} \right|_{\mathrm{Inj\, Mult = 1}}\\
& = 0.5724 \cdot \sqrt{4} = 1.1448 \mathrm{\, per\, unit}\end{align}

The voltage margin at a power injection multiplier of four is $\sigma = \frac{1.0350}{1.1448} = 0.9041$. $\sigma < 1$ indicates that there is no solution at a power injection multiplier of four. To potentially achieve a power flow solution, the slack bus voltage must increase by at least a factor of $\frac{1}{0.1098} = 1.1061$ (with corresponding proportional increases in all PV bus voltages). The grey PV curve in Fig. \ref{f:pv_hv_118} has the voltages thus increased. Since no solutions are evident from the PV curve at an injection multiplier of four, it appears that this is not a large enough voltage increase to obtain solvability. This is a result of the fact that we use a sufficient condition for power flow insolvability to calculate the voltage margin; failing to satisfy the the sufficient condition for power flow insolvability does not ensure the existence of a solution.

The power injection margin $\eta$ can also be calculated at a power injection multiplier of four using \eqref{pinj_margin_def}.

\begin{align}\nonumber \eta & =  \left( \frac{V_{0}}{\left. \underline{V}_{slack}^{min} \right|_{\mathrm{Inj\, Mult = 4}}} \right)^2 \\
& = \left(\frac{1.0350}{1.1448} \right)^2 = 0.8174 \end{align}

$\eta < 1$ implies that no solution exists at a power injection multiplier of four. The power injection multiplier also indicates that no solution can exist for power injection multipliers greater than $ 0.8174 \cdot 4 = 3.2695$. The ``nose'' point of the black (nominal) PV curve in Fig. \ref{f:pv_hv_14} is slightly lower than this upper limit on power flow solvability.

%
%
%
%
%
%

\section{Solutions with Non-Zero Duality Gap}
\label{l:nonzero_gap}

If the $\mathbf{A}$ matrix in \eqref{cond1dualA} has nullspace with rank greater than two, a non-zero duality gap exists between the primal problem \eqref{cond1} and the dual problem \eqref{cond1dual} \cite{lavaei}. The dual problem will give a $\underline{V}_{slack}^{min}$ that is less than the $V_{slack}^{min}$ of the primal problem. If $\underline{V}_{slack}^{min} \leq  V_0 < V_{slack}^{min}$, condition \eqref{cond1leq} may hold without the existence of a power flow solution. Since directly calculating $V_{slack}^{min}$ is difficult due to the non-convexity of the power flow equations in \eqref{cond1}, it is hard to determine whether a solution exists for cases where \eqref{cond1leq} holds and $\mathrm{rank}\left(\mathrm{null}\left(\mathbf{A}\right)\right) > 2$. The IEEE 118-bus system in Section \ref{l:118bus_cond1} provides an example of such a case.

As the controlled voltages are decreased, the solutions to the power flow equations come together in bifurcations. We conjecture that the rank of the nullspace of the $\mathbf{A}$ matrix depends on where these solutions bifurcate. The minimum slack bus voltage of the primal optimization \eqref{cond1} is obtained when the last solutions bifurcate. If the minimum slack bus voltage is obtained when two solutions bifurcate far from the point of bifurcation of any previously existing solutions, we conjecture that $\mathrm{rank}\left(\mathrm{null}\left(\mathbf{A}\right)\right) \leq 2$. If more than two solutions bifurcate at or near the same minimum slack bus voltage (e.g., the final two pairs of solutions bifurcate with the same or similar values of $V_{slack}^{min}$), we conjecture that $\mathrm{rank}\left(\mathrm{null}\left(\mathbf{A}\right)\right) > 2$. 

\begin{figure}[b]
\centering
\includegraphics[totalheight=0.08\textheight]{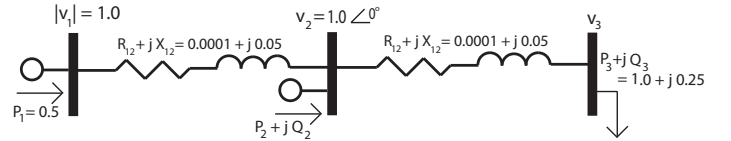}
\caption{Three-Bus System}
\label{f:threebussystemdiagram}
\end{figure}

The system diagram in Fig. \ref{f:threebussystemdiagram} provides an example with two pairs of solutions bifurcating at the same value of minimum slack bus voltage, resulting in $\mathrm{rank}\left(\mathrm{rank}\left( \mathbf{A}\right) \right) = 4$. Bus one is a PV bus with voltage magnitude $1.0$ per unit and generation of 0.5 per unit, bus two is a slack bus with voltage $1.0 \angle 0^\circ$, and bus three is a PQ bus with load $1.0 + j 0.25$ per unit. 

\begin{figure}[t]
\centering
\includegraphics[totalheight=0.28\textheight]{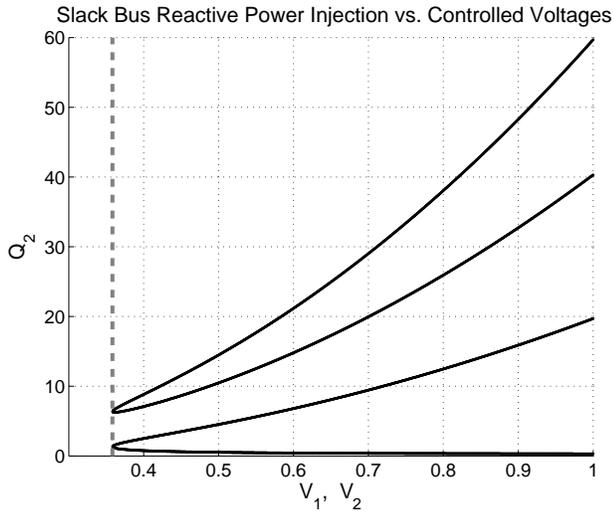}
\caption{Solutions for Three-Bus System}
\label{f:threebusbifurcation}
\end{figure}

Fig. \ref{f:threebusbifurcation} shows the slack bus reactive power injection $Q_2$ for all four existing solutions as the controlled voltages at buses one and two are uniformly decreased. When the controlled voltage magnitudes reach 0.3582 per unit, both pairs of solutions bifurcate. Note that both pairs of solutions bifurcate at different values of slack bus reactive power injection. That is, even though both solution pairs disappear at the same value of controlled voltages, all four solutions do not come together to the same point.


For this example, $\underline{V}_{slack}^{min}$ from optimization \eqref{cond1dual} finds the correct minimum slack bus voltage of 0.3582 per unit. The bifurcation points for both solution pairs in Fig. \ref{f:threebusbifurcation} have this voltage. Thus, the optimization problem \eqref{cond1dual} has no bias toward choosing either of these two solutions. Indeed, since all solutions must be connected in the convex dual formulation, any matrices in the dual feasible space that lie between these solutions, including matrices that have nullspaces with high ranks, are possible solutions to \eqref{cond1dual}. Thus, the optimization problem \eqref{cond1dual} is not expected to obtain $\mathrm{rank}\left(\mathrm{null}\left( \mathbf{A}\right)\right) \leq 2$. We conjecture that similar explanations apply to other systems with $\mathrm{rank}\left(\mathrm{null}\left( \mathbf{A}\right)\right) > 2$, such as the IEEE 118-bus system used in Section \ref{l:118bus_cond1}.

Using optimization codes YALMIP \cite{yalmip} and SeDuMi \cite{sedumi}, we obtained a solution to \eqref{cond1dual} for the three-bus system with $\mathrm{rank}\left(\mathrm{null}\left( \mathbf{A}\right)\right) = 4$. For this system, adding a small bias in the objective function to favor larger or smaller slack bus reactive power injection resulted in solutions with $\mathrm{rank}\left(\mathrm{null}\left( \mathbf{A}\right)\right) = 2$. Appropriate biases to obtain solutions with $\mathrm{rank}\left(\mathrm{null}\left( \mathbf{A}\right)\right) = 2$ have not yet been found for other systems (e.g., the IEEE 118-bus system used in Section \ref{l:118bus_cond1}).

The sufficient condition \eqref{cond1g} does not generally give a solution to the power flow equations with slack bus voltage equal to the specified value of $V_0$; traditional techniques are required to actually find a potentially existing solution. In an attempt to actually find power flow solutions, a modification to the objective function of \eqref{cond1} was attempted. Specifically, the objective function \eqref{cond1obj} was replaced with \eqref{cond2obj},

\begin{equation}\label{cond2obj} \mathrm{min}\quad \left(V_{slack}^2 - V_{0}^2 \right)^2 \end{equation}

This modified objective function minimizes the squared difference in squared slack bus voltage magnitude (squared voltage magnitudes are necessary to easily form the dual problem). The power flow equations have no solution when the optimal objective value is greater than zero.

Since the constraints \eqref{cond1P}, \eqref{cond1Q}, \eqref{cond1V} are unchanged, the existence proof given in Section \ref{l:solution_existence_proof} applies, so this modified condition can be evaluated. 

In the same manner as the original sufficient condition, the dual of the optimization problem with the modified objective function \eqref{cond2obj} provides a lower bound on the primal objective function. If the dual objective value is greater than zero, no solutions exist for the power flow equations. If, additionally, the $\mathbf{A}$ matrix  has a nullspace with rank less than or equal to two, the solution to the power flow equations can be obtained \cite{lavaei}.

Unfortunately, the optimization problem with this modified objective function generally failed to provide $\mathbf{A}$ matrices that had nullspaces with appropriate rank on all IEEE power flow test systems \cite{ieee_test_cases}. Also, the modified formulation often had numeric problems. This modified objective function \eqref{cond2obj} was thus not practically useful. Future work includes investigating other modifications that allow for determining the solution to the power flow equations in addition to providing a sufficient condition of solution non-existence.

\section{Conclusion and Future Work}
\label{l:conclusion}

We have presented a sufficient condition for identifying insolvability of the power flow equations. This sufficient condition requires the evaluation of an optimization problem. We have proven that this optimization problem is feasible for lossless power systems and argued that practical power systems should also yield a feasible optimization problem. In order to quantify the degree of solvability, we developed controlled voltage and power injection margins from the sufficient condition that provide upper bounds on the distance to the power flow solvability boundary. Finally, we applied the sufficient conditions, voltage margin, and power injection margin to the IEEE 14-bus and 118-bus example systems. The IEEE 118-bus system exemplified the fact that this is a \emph{sufficient} condition for power flow insolvability; it is possible to find cases that appear to have no solution even though the sufficient condition for insolvability is not satisfied. However, the majority of power systems we investigated yielded results similar to the IEEE 14-bus system where a power flow solution was found with a Newton-Raphson algorithm up to the point identified by the sufficient condition as insolvable. 

There are several open questions associated with this work that we intend to pursue. We will investigate how the voltage and power injection margins change with non-uniform power injection scaling. For instance, what happens if the power injections are scaled at a subset of the buses, or if the power injections are scaled with non-constant power factor? The sensitivity of $V_{slack}^{min}$ in \eqref{cond1} to the active and reactive power flow constraints \eqref{cond1P} and \eqref{cond1Q} (as determined by the Lagrange multipliers $\lambda$ and $\gamma$, respectively) provides insight for small changes, but numerical experience indicates that these sensitivities are only valid for very small perturbations. The impact of large power injection changes requires further investigation.


Another area open to further study regards computational issues. The size of positive semidefinite constraint in the dual optimization problem \eqref{cond1dualA} scales as $\left(2n\right)^2$, where $n$ is the number of buses in the system. This scaling prohibits directly solving \eqref{cond1dual} for large systems. Future work includes researching whether the structure and sparsity of the power flow equations can be exploited to solve this optimization problem quickly and efficiently for large systems.

Finally, we intend to further investigate our conjectured explanation of solutions with non-zero duality gap. Do multiple pairs of solutions bifurcate near the minimum slack bus voltage for other systems? What system topologies yield the possibility of having $\mathrm{rank}\left(\mathrm{null}\left( \mathbf{A}\right)\right) > 2$?

\section*{Acknowledgment}
The authors acknowledge support of this work by U.S. Department of Energy under award \#DE-SC0002319, 
as well as by the National Science Foundation under IUCRC award \#0968833. Daniel Molzahn acknowledges the support of the National Science Foundation Graduate Research Fellowship.



\IEEEtriggeratref{12}


\bibliographystyle{IEEEtran}
\bibliography{IEEEabrv,pfconditionrefs}

\begin{thebibliography}{10}
\providecommand{\url}[1]{#1}
\csname url@samestyle\endcsname
\providecommand{\newblock}{\relax}
\providecommand{\bibinfo}[2]{#2}
\providecommand{\BIBentrySTDinterwordspacing}{\spaceskip=0pt\relax}
\providecommand{\BIBentryALTinterwordstretchfactor}{4}
\providecommand{\BIBentryALTinterwordspacing}{\spaceskip=\fontdimen2\font plus
\BIBentryALTinterwordstretchfactor\fontdimen3\font minus
  \fontdimen4\font\relax}
\providecommand{\BIBforeignlanguage}[2]{{%
\expandafter\ifx\csname l@#1\endcsname\relax
\typeout{** WARNING: IEEEtran.bst: No hyphenation pattern has been}%
\typeout{** loaded for the language `#1'. Using the pattern for}%
\typeout{** the default language instead.}%
\else
\language=\csname l@#1\endcsname
\fi
#2}}
\providecommand{\BIBdecl}{\relax}
\BIBdecl

\bibitem{saadat}
H.~Saadat, \emph{{Power System Analysis}}.\hskip 1em plus 0.5em minus
  0.4em\relax McGraw-Hill, 2005.

\bibitem{glover_sarma_overbye}
J.~Glover, M.~Sarma, and T.~Overbye, \emph{{Power System Analysis and
  Design}}.\hskip 1em plus 0.5em minus 0.4em\relax Thompson Learning, 2008.

\bibitem{baillieul_brynes}
J.~Baillieul and C.~Byrnes, ``{Geometric Critical Point Analysis of Lossless
  Power System Models},'' \emph{{IEEE Transactions on Circuits and Systems}},
  vol.~29, no.~11, pp. 724--737, Nov 1982.

\bibitem{thorp1989}
J.~Thorp and S.~Naqavi, ``Load flow fractals,'' in \emph{Proceedings of the
  28th IEEE Conference on Decision and Control, 1989.}, dec 1989, pp. 1822
  --1827 vol.2.

\bibitem{garcia_zangwill}
W.~Zangwill and C.~Garcia, \emph{{Pathways to Solutions, Fixed Points, and
  Equilibria}}.\hskip 1em plus 0.5em minus 0.4em\relax Prentice-Hall, 1981.

\bibitem{DOE1}
J.~Jarjis and F.~D. Galiana, ``{Analysis and Characterization of Security
  Regions in Power Systems, Part I: Load Flow Feasibility Conditions in Power
  Systems},'' McGill University, Final Report, U.S. Department of Energy,
  Division of Electric Energy Systems, DOE/ET/29108--T1-Pt.1, March 1980.

\bibitem{ilic1992}
M.~Ilic, ``{Network Theoretic Conditions for Existence and Uniqueness of Steady
  State Solutions to Electric Power Circuits},'' in \emph{Proceedings of 1992
  IEEE International Symposium on Circuits and Systems (ISCAS)}, vol.~6, May
  1992, pp. 2821--2828 vol.6.

\bibitem{test_illcond}
M.~Ebrahimpour and J.~Dorsey, ``{A Test for the Existence of Solutions in
  Ill-Conditioned Power Systems},'' in \emph{IEEE Proceedings of Southeastcon
  '91.}, vol.~1, Apr 1991, pp. 444--448.

\bibitem{necessary_stsl}
S.~Grijalva, ``{Individual Branch and Path Necessary Conditions for Saddle-Node
  Bifurcation Voltage Collapse},'' \emph{{IEEE Transactions on Power Systems}},
  vol.~27, no.~1, pp. 12--19, Feb. 2012.

\bibitem{hiskens2001}
I.~Hiskens and R.~Davy, ``{Exploring the Power Flow Solution Space Boundary},''
  \emph{{IEEE Transactions on Power Systems}}, vol.~16, no.~3, pp. 389--395,
  Aug 2001.

\bibitem{optimal_multiplier}
S.~Iwamoto and Y.~Tamura, ``{A Load Flow Calculation Method for Ill-Conditioned
  Power Systems},'' \emph{{IEEE Transactions on Power Apparatus and Systems}},
  vol. PAS-100, no.~4, pp. 1736--1743, April 1981.

\bibitem{overbye1994}
T.~Overbye, ``{A Power Flow Measure for Unsolvable Cases},'' \emph{{IEEE
  Transactions on Power Systems}}, vol.~9, no.~3, pp. 1359--1365, Aug 1994.

\bibitem{overbye1995}
------, ``{Computation of a Practical Method to Restore Power Flow
  Solvability},'' \emph{{IEEE Transactions on Power Systems}}, vol.~10, no.~1,
  pp. 280--287, Feb 1995.

\bibitem{dobson1993}
I.~Dobson and L.~Lu, ``{New Methods for Computing a Closest Saddle Node
  Bifurcation and Worst Case Load Power Margin for Voltage Collapse},''
  \emph{{IEEE Transactions on Power Systems}}, vol.~8, no.~3, pp. 905--913, Aug
  1993.

\bibitem{alvarado1994}
F.~Alvarado, I.~Dobson, and Y.~Hu, ``{Computation of Closest Bifurcations in
  Power Systems},'' \emph{{IEEE Transactions on Power Systems}}, vol.~9, no.~2,
  pp. 918--928, May 1994.

\bibitem{continuation_optimization_pf}
F.~Echavarren, E.~Lobato, and L.~Rouco, ``{A Power Flow Solvability
  Identification and Calculation Algorithm},'' \emph{Electric Power Systems
  Research}, vol.~76, no.~4, pp. 242--250, 2006.

\bibitem{interior_point_unsolvable}
S.~Granville, J.~Mello, and A.~Melo, ``{Application of Interior Point Methods
  to Power Flow Unsolvability},'' \emph{{IEEE Transactions on Power Systems}},
  vol.~11, no.~2, pp. 1096--1103, May 1996.

\bibitem{lesieutre2008}
V.~Donde, V.~Lopez, B.~Lesieutre, A.~Pinar, C.~Yang, and J.~Meza, ``{Severe
  Multiple Contingency Screening in Electric Power Systems},'' \emph{{IEEE
  Transactions on Power Systems}}, vol.~23, no.~2, pp. 406--417, {May} 2008.

\bibitem{feng1998}
Z.~Feng, V.~Ajjarapu, and D.~Maratukulam, ``{A Practical Minimum Load Shedding
  Strategy to Mitigate Voltage Collapse},'' \emph{{IEEE Transactions on Power
  Systems}}, vol.~13, no.~4, pp. 1285--1290, Nov 1998.

\bibitem{pf_insolvability_comparison}
J.~Zhao, Y.~Wang, and P.~Ju, ``{Evaluation of Methods for Measuring the
  Insolvability of Power Flow},'' in \emph{{Third International Conference on
  Electric Utility Deregulation and Restructuring and Power Technologies, 2008
  (DRPT 2008)}}, April 2008, pp. 920--925.

\bibitem{sdp_china}
X.~Bai, H.~Wei, K.~Fujisawa, and Y.~Wang, ``{Semidefinite Programming for
  Optimal Power Flow Problems},'' \emph{International Journal of Electrical
  Power \& Energy Systems}, vol.~30, no. 6-7, pp. 383--392, 2008.

\bibitem{lavaei}
J.~Lavaei and S.~Low, ``{Zero Duality Gap in Optimal Power Flow Problem},''
  \emph{{IEEE Transactions on Power Systems}}, vol.~27, no.~1, pp. 92--107,
  Feb. 2012.

\bibitem{allerton2011}
B.~C. Lesieutre, D.~K. Molzahn, A.~R. Borden, and C.~L. DeMarco, ``{Examining
  the Limits of the Application of Semidefinite Programming to Power Flow
  Problems},'' in \emph{{49th Annual Allerton Conference on Communication,
  Control, and Computing, 2011}}, Sept. 28-30 2011.

\bibitem{bernie_opfconvexity}
B.~Lesieutre and I.~Hiskens, ``{Convexity of the Set of Feasible Injections and
  Revenue Adequacy in FTR Markets},'' \emph{IEEE Transactions on Power
  Systems}, vol.~20, no.~4, pp. 1790 -- 1798, November 2005.

\bibitem{Chua_Desoer_Kuh}
L.~Chua, C.~Desoer, and E.~Kuh, \emph{{Linear and Nonlinear Circuits}}.\hskip
  1em plus 0.5em minus 0.4em\relax McGraw-Hill, New York, 1987.

\bibitem{ieee_test_cases}
\BIBentryALTinterwordspacing
{Power Systems Test Case Archive}. {University of Washington Department of
  Electrical Engineering}. [Online]. Available:
  \url{http://www.ee.washington.edu/research/pstca/}
\BIBentrySTDinterwordspacing

\bibitem{hirsch_smale}
M.~Hirsch and S.~Smale, \emph{{Differential Equations, Dynamical Systems, and
  Linear Algebra}}.\hskip 1em plus 0.5em minus 0.4em\relax Academic Press New
  York, 1974.

\bibitem{horn1985}
R.~Horn and C.~Johnson, \emph{{Matrix Analysis}}.\hskip 1em plus 0.5em minus
  0.4em\relax Cambridge University Press, 1985.

\bibitem{yalmip}
J.~Lofberg, ``{YALMIP: A Toolbox for Modeling and Optimization in MATLAB},'' in
  \emph{{IEEE International Symposium on Computer Aided Control Systems Design,
  2004}}.\hskip 1em plus 0.5em minus 0.4em\relax IEEE, 2004, pp. 284--289.

\bibitem{sedumi}
J.~Sturm, ``{Using SeDuMi 1.02, A MATLAB Toolbox for Optimization Over
  Symmetric Cones},'' \emph{{Optimization Methods and Software}}, vol.~11,
  no.~1, pp. 625--653, 1999.

\bibitem{cpf}
V.~Ajjarapu and C.~Christy, ``{The Continuation Power Flow: A Tool for Steady
  State Voltage Stability Analysis},'' \emph{{IEEE Transactions on Power
  Systems}}, vol.~7, no.~1, pp. 416--423, Feb 1992.

\end{thebibliography}
\end{document}